\documentclass[a4paper,11pt]{article}
\usepackage{amsmath,amssymb}
\usepackage{latexsym}
\usepackage{graphicx}

\usepackage{amssymb,amsfonts,amsmath,color,
latexsym, epsfig,cite, psfrag,eepic,colordvi}

\newtheorem{theorem}{Theorem}
\newtheorem{lemma}[theorem]{Lemma}

\newtheorem{definition}[theorem]{Definition}
\newtheorem{problem}[theorem]{Problem}

\newcommand{\qed}{~$\Box $}

\begin{document}

\title{Note\\
On the Complexity of All $(g,f)$-Factors Problem}

\author{Hongliang Lu, Wei Wang and Yaolin Jiang
\\ School of Mathematics and Statistics\\
Xi'an Jiaotong University\\
Xi'an, Shaanxi 710049, China\\
}

\maketitle

\begin{abstract}
Let $G$ be a graph with vertex set $V$ and let $g, f : V\rightarrow \mathbb{Z}^+$ be two functions such that $g\le f$. We say that $G$ has all
$(g, f )$-factors if $G$ has an $h$-factor for every $h: V\rightarrow \mathbb{Z}^+$ such that $g(v)\le h(v)\le f (v)$
for every $v\in V$ and $\sum_{v\in{V}}h(v)\equiv 0\pmod 2$. Two decades ago, Niessen derived from
Tutte's $f$-factor theorem a similar characterization for the property of graphs having
all $(g, f )$-factors and asked whether there is a
polynomial time algorithm for testing whether a graph $G$ has all $(g, f )$-factors (A characterization of graphs having all $(g, f )$-Factors, \emph{J. Combin. Theory, Ser. B},  \textbf{72} (1998), 152--156). In this paper, we show that it is NP-hard to determine whether a graph $G$ has all $(g,f)$-factors, which gives a negative answer to the question of Niessen.
\end{abstract}

\section{Introduction}
We consider finite simple graphs, i.e., undirected graph without loops or multiple edges. Let $G$ be a graph with vertex set $V(G)$ and edge set $E(G)$.

Let $\omega(G)$ denote the number of components of a graph $G$. A graph $G$ is \emph{$t$-tough}
if $|S| \le t\omega(G-S)$ for every subset $S$ of the vertex set $V (G)$ with $\omega(G-S) > 1$. The
\emph{toughness} of $G$, denoted $t(G)$, is the maximum value of $t$ for which $G$ is $t$-tough (taking
$t(K_n)=\infty$ for all $n \le 1$).
We say that a graph $G$ is \emph{almost $1$-tough} if for any $S\subseteq V(G)$,
\begin{align*}
\omega(G-S)\le |S|+1.
\end{align*}

Let $H:V(G)\rightarrow 2^{\mathbb{Z}^+}$ be a set function. An $H$-factor $F$ is a spanning subgraph such that $d_F(v)\in H(v)$ for every $v\in V(G)$.
Given two integer-valued functions $g,f:V(G)\rightarrow \mathbb{Z}^+$ such that $g(v)\le f(v)$ for all $v\in V(G)$, an $H$-factor is also called \emph{$(g,f)$-factor} if $H(v)=[g(v),f(v)]$ for all $v\in V(G)$. In particular, if $g= f$, then an $(f,f)$-factor is also called an \emph{$f$-factor.}
We say that $G$ has all $(g,f)$-factors if for any integer-valued function defined on $V(G)$ such that $h(V(G))\equiv 0\pmod 2$ and $g(v)\le h(v)\le f(v)$
for all $v\in V(G)$, $G$ contains an $h$-factor.  We write
\[
 H^{-1}(1):=\{v\in V(G): H(v)=\{1\}\}.
\]
For simplicity, we write $g^{-1}(1):=H^{-1}(1)$ if $H(v)=\{g(v)\}$ for all $v\in V(G)$.

Niessen~\cite{Ni98} derived from
Tutte's $f$-factor theorem a similar characterization for the property of graphs having
all $(g, f )$-factors.
\begin{theorem}[Niessen~\cite{Ni98}]
$G$ has all $(g, f )$-factors if and only if
\begin{equation*}
  g(D)-f (S)+d_{G-D}(S)-q_G(D, S, g, f )\ge\left\{
  \begin{array}{ll}
    -1, & \hbox{$g\neq f$;} \\
    0, & \hbox{$g=f$.}
  \end{array}
\right.
\end{equation*}
for all disjoint sets $D, S\subseteq V$, where $q_G(D, S, g, f )$ denotes the number of
components $C$ of $G-(D\cup S)$ such that there exists a vertex $v \in V(C)$ with
$g(v)< f (v)$ or $e_G(V(C), S)+ f (V(C))\equiv 1 \pmod 2$.
\end{theorem}

It is well-known that there exists a polynomial time algorithm to determine whether a graph has a $(g,f)$-factor.
An open problem was naturally proposed in \cite{Ni98}:
\begin{problem}\label{prob1}
 Is there a polynomial time algorithm for testing whether a graph $G$ has all $(g, f )$-factors?
\end{problem}

In this note, we obtain the following theorem, which gives a negative answer to Problem \ref{prob1}, unless P=NP.

\begin{theorem}\label{main-thm}
It is NP-hard to determine whether a graph $G$ has all $(g, f )$-factors.
\end{theorem}

The main ingredient in the proof of Theorem~\ref{main-thm} is the following two results.
\begin{theorem}[Bauer et. al.,~\cite{BHMS1}]\label{NPC-1t}
It is NP-hard to recognize
1-tough cubic graphs.
\end{theorem}

\begin{theorem}[Kano and Lu,~\cite{KaLu}] \label{KaLu}
Let $G$ be a connected graph. A graph $G$ has an $H$-factor for every $H:V(G)\rightarrow \{\{1\},\{0,2\}\}$  with $|H^{-1}(1)|$ even if and only if
$G$ is almost 1-tough.
\end{theorem}

Theorem \ref{KaLu} gives a polynomial time reduction from the toughness of graphs to the degree constrained factors.

\section{The Proof of Theorem~\ref{main-thm}}

To prove Theorem~\ref{main-thm}, we need the following notations.
\begin{definition}\label{N1}
 Let $G$ be a graph with $V(G)=\{v_1,\ldots,v_n\}$. Write $X=\{x_1,\ldots,x_n\}$, $Y=\{y_1,\ldots,y_n\}$.
 \begin{itemize}
\item [(i)] For any vertex $x$ of $G$, let $G^x$ denote the graph obtained from $G$
by adding a new vertex $x'$ together with a new edge $xx'$, that is, $G^x=G+xx'$.

\item [(ii)] Let $G_L$ be a graph  with  vertex set $V(G_L)$ and edge set $E(G_L)$, where $V(G_L)=V(G)\cup X\cup Y$ and $E(G_L)=E(G)\cup \{x_iy_i,y_iv_i, v_ix_i\ |\ i\in [n]\}$;

 \item [(iii)] For any $H: V(G)\rightarrow \{\{1\},\{0,2\}\}$,  we define $h_H: V(G_L)\rightarrow \{1,2\}$ such that
\begin{equation*}
  h_H(u)=\left\{
     \begin{array}{ll}
       2, & \hbox{if $H(u)=\{0,2\}$ and $u\in V(G)$;} \\
       1, & \hbox{otherwise.}
     \end{array}
   \right.
\end{equation*}
 \end{itemize}
\end{definition}

Now we show the following result.

\begin{lemma}\label{t-at}
A graph $G$ is 1-tough if and only if for any $x\in V(G)$, $G^x$ is almost 1-tough.
\end{lemma}

\noindent\textbf{Proof.} Necessity. Suppose that $G$ is 1-tough. Let $x\in V(G)$. For any $S\subseteq V(G^x)$,
\begin{align*}
\omega(G^x-S)\le \omega(G-S)+|\{x'\}|\le |S|+1.
\end{align*}
Sufficiency. By contradiction, suppose that $G$ is not 1-tough. Then there exists $\emptyset\neq S\subseteq V(G)$ such that $\omega(G-S)\ge |S|+1$. Let $x\in S$. One can see that
\begin{align*}
\omega(G^x-S)=\omega(G-S)+1\ge |S|+2,
\end{align*}
which contradicts the fact that $G$ is 1-tough. \qed

From Theorem \ref{NPC-1t} and Lemma \ref{t-at}, one may see that
\begin{lemma}\label{Coro1}
It is NP-hard to recognize almost 1-tough cubic graph.
\end{lemma}

\begin{lemma}\label{lem2}
Let $G$ be a connected cubic graph and let $H:V(G)\rightarrow \{\{1\},\{0,2\}\}$. Then $G$ contains an $H$-factor if and only if $G_L$ contains an $h_H$-factor.
\end{lemma}

\noindent\textbf{Proof.} Necessity. Suppose that $G$ contains an $H$-factor $F$.
Define $M_1=\{x_iy_i\ |\ d_F(v_i)\in \{1, 2\}\ \mbox{for $i\in [n]$}\}$ and $M_2=\{x_iv_i,y_iv_i\ |\ d_F(v_i)=0\ \mbox{for $i\in [n]$}\}$.
Let $F'$ be  a spanning subgraph of $G_L$ with edge set $E(F)\cup M_1\cup M_2$. From the definition of function $h_H$, one can see that
\begin{equation*}
  h_H(x)=d_{F'}(x)=\left\{
     \begin{array}{ll}
       2, & \hbox{if $H(x)=\{0,2\}$ and $x\in V(G)$;} \\
       1, & \hbox{otherwise.}
     \end{array}
   \right.
\end{equation*}
So $F'$ is an $h_H$-factor of $G_L$.

Sufficiency. Suppose that $F'$ be an $h_H$-factor of $G_L$. Let $F$ be a spanning subgraph of $G$ with edge set $E(F')\cap E(G)$. Now we show that $F$ is an $H$-factor of $G$. From Definition \ref{N1} (iii), one can see that
\begin{equation*}
  H(u)=\left\{
     \begin{array}{ll}
       \{1\}, & \hbox{$h_H(u)=1$ and $u\in V(G)$;} \\
       \{0,2\}, & \hbox{$h_H(u)=2$ and $u\in V(G)$.}
     \end{array}
   \right.
\end{equation*}
Consider $h_H(v_i)=d_{F'}(v_i)=1$. Since $d_{F'}(x_i)=d_{F'}(y_i)=1$, then we have $x_iv_i,y_iv_i\notin E(F')$. So we have $d_F(x_i)=1\in H(v_i)$. Next we may assume that $d_{F'}(v_i)=2$. Then either $\{x_iv_i,y_iv_i\}\subseteq E(F')$ or $\{x_iv_i,y_iv_i\}\cap E(F')=\emptyset$. In the former case, we have $d_F(v_i)=0$, and in the latter case we have $d_{F}(v_i)=2$.  So in both cases, one can see that $d_F(v_i)\in \{0,2\}$ and $F$ is an $H$-factor of $G$.
This completes the proof. \qed

\begin{lemma}\label{lem3}
Let $G$ be a connected cubic graph. Then $G$ contains an $H$-factor for any $H:V(G)\rightarrow \{\{1\},\{0,2\}\}$ with $|H^{-1}(1)|$ even if and only if $G_L$ contains all $(g,f)$-factors, where
$g\equiv 1$ and \begin{equation*}
  f(u)=\left\{
     \begin{array}{ll}
       2, & \hbox{if $u\in V(G)$;} \\
       1, & \hbox{else $u\in V(G_L)-V(G)$.}
     \end{array}
   \right.
\end{equation*}

\end{lemma}

\noindent\textbf{Proof.} Define
\begin{equation*}
\mathcal{H}=\{H\ |\ \mbox{$H: V(G)\rightarrow  \{\{1\},\{0,2\}\}$ such that $|H^{-1}(1)|$ is even}\}.
\end{equation*}
and
\begin{equation*}
\mathcal{F}=\{h\ |\ h: V(G_L)\rightarrow  \{1,2\} \\ \mbox{ such that $|h^{-1}(1)|$ is even
and $h(v)=1$ $\forall$ $v\in X\cup Y$}\}.
\end{equation*}
Let $J:\mathcal{H}\rightarrow \mathcal{F}$ such that $J(H)=h_H$ for all $H\in \mathcal{H}$. From the definition of $h_H$, one can see that $J$ is well-defined.
By Lemma \ref{lem2}, it is suffices for us to show that $J$
  is a bijection.  Firstly, we show that $J$ is injective. For any $H,H'\in \mathcal{H}$ such that  $H\neq H'$, there exists $x\in V(G)$, such that $H(x)\neq H'(x)$. Without loss of  generality, we may assume that $H(x)=\{1\}$ and $H'(x)=\{0,2\}$. By Definition \ref{N1} (iii), we have $h_H(x)=1$ and $h_{H'}(x)=2$, which implies $h_H\neq h_{H'}$.
Next we show $J$ is a surjection. For every $h\in \mathcal{F}$, we may define $H:V(G)\rightarrow \{\{1\},\{0,2\}\}$ such that
\begin{equation*}
  H(v)=\left\{
     \begin{array}{ll}
       \{1\}, & \hbox{if $h(v)=1$ and $v\in V(G)$;} \\
       \{0,2\}, & \hbox{else $h(v)=2$ and $v\in V(G)$.}
     \end{array}
   \right.
\end{equation*}
Since $|h^{-1}(1)|$ is even, $|H^{-1}(1)|=|h^{-1}(1)\cap V(G)|$ is even. So we get that $H\in \mathcal{H}$ and $J(H)=h$. Thus $J$ is surjective.

 This completes the proof. \qed\\

\noindent
\textbf{Proof of Theorem \ref{main-thm}.}
Let $G$ be a connected cubic graph.  By Theorem~\ref{KaLu}, $G$ is almost 1-tough  if and only if
$G$ contains an $H$-factor for  any $H:V(G)\rightarrow \{\{1\},\{0,2\}\}$ with $|H^{-1}(1)|$ even.
Thus by Lemma \ref{lem3}, one can see that $G$ is almost 1-tough if and only if
$G_L$ contains all $(g,f)$-factors, where
$g\equiv 1$ and \begin{equation*}
  f(u)=\left\{
     \begin{array}{ll}
       2, & \hbox{if $u\in V(G)$;} \\
       1, & \hbox{else $u\in V(G_L)-V(G)$.}
     \end{array}
   \right.
\end{equation*}
Since it is NP-hard to recognize almost 1-tough graphs by Lemma~\ref{Coro1}, one can see that it is NP-hard to determine whether a graph contains all $(g,f)$-factors. This completes the proof of Theorem~\ref{main-thm}. \qed

\end{document}